\documentclass[12pt]{article}

\usepackage{theorem,amssymb}

\advance \topmargin by -\headheight
\advance \topmargin by -\headsep
\evensidemargin \oddsidemargin
\author{J.-P. Allouche\thanks{partially supported 
by MENESR, ACI NIM 154 Num\'eration.} \\
CNRS, LRI, B\^atiment 490 \\
F-91405 Orsay Cedex \\
France\\
{\tt allouche@lri.fr} \\
\and
J. Shallit \\
School of Computer Science, University of Waterloo \\
Waterloo, Ontario  N2L 3G1 \\
Canada \\
{\tt shallit@graceland.uwaterloo.ca} \\
\and
J. Sondow \\
209 West 97th Street \\
New York, NY 10025  \\
USA\\
{\tt jsondow@alumni.princeton.edu} \\
}

\title{Summation of Series Defined by Counting Blocks of Digits}

\date{ }

\def \proof{\bigbreak\noindent{\it Proof.\ \ }}

\def \endpf{{\ \ $\Box$ \medbreak}}

\newtheorem{theorem}{Theorem}
\newtheorem{lemma}{Lemma}

\theorembodyfont{\rm}

\newtheorem{remark}{Remark}
\newtheorem{example}{Example}

\begin{document}

\maketitle

\begin{abstract}
We discuss the summation of certain series defined by counting blocks of 
digits in the $B$-ary expansion of an integer. For example, if $s_2(n)$ 
denotes the sum of the base-$2$ digits of $n$, we show that 
$\sum_{n \geq 1} s_2(n)/(2n(2n+1)) = (\gamma + \log \frac{4}{\pi})/2$.
We recover this previous result of Sondow and provide several generalizations.

\medskip

\end{abstract}

\section{Introduction} 

A classical series with rational terms, known as Vacca's series
\cite{Vacca} or in an equivalent integral form as Catalan's integral
\cite{Catalan} (see also \cite{BerBow} and \cite{SondZudi}), evaluates 
to Euler's constant $\gamma$:
$$
\gamma = \sum_{n \geq 1} \frac{(-1)^n}{n} 
\left\lfloor \frac{\log n}{\log 2} \right\rfloor = 
\int_0^1 \frac{1}{1+x} \sum_{n \geq 1} x^{2^n -1}dx.
$$
In a recent paper \cite{Sondow2} Sondow gave the following two 
formulas:
$$
\gamma^{\pm} = \sum_{n \geq 1} \frac{N_1(n) \pm N_0(n)}{2n(2n+1)}
$$
where $\gamma^+ = \gamma$ is the Euler constant, 
$\gamma^- = \log \frac{4}{\pi}$ is the ``alternating Euler 
constant'' \cite{Sondow1}, and $N_1(n)$ (resp. $N_0(n)$) is 
the number of $1$'s (resp. $0$'s) in the binary expansion of the 
integer $n$. The series for $\gamma^+ = \gamma$ is equivalent to 
Vacca's. The formulas for $\gamma^{\pm}$ show in particular that
$$
\sum_{n \geq 1} \frac{s_2(n)}{2n(2n+1)} =
\frac{\gamma + \log \frac{4}{\pi}}{2}
$$
where $s_2(n)$ is the sum of the binary digits of the integer $n$.

This last formula reminds us of one of the problems posed at the 1981 Putnam 
competition \cite{KAH}: Determine whether or not 
$$
\exp\left(\sum_{n \geq 1} \frac{s_2(n)}{n(n+1)}\right)
$$
is a rational number. 
In fact, $\sum \frac{s_2(n)}{n(n+1)} = 2\log 2$. A generalization 
was proven by Shallit \cite{Shallit}, where the base $2$ is replaced by
any integer base $B \geq 2$. A more general result, where the sum of digits is
replaced by the function
$N_{w,B}(n)$, which counts the number of occurrences of the
block $w$ in the $B$-ary expansion of the integer $n$, was given 
by Allouche and Shallit \cite{AS}.

The purpose of the present paper is to show that the result of \cite{Sondow2} 
cited above can be deduced from a general lemma in \cite{AS}. Furthermore, we 
sum the series
$$
\sum_{n \geq 1} \frac{N_{w,2}(n)}{2n(2n+1)}
\ \ \ \mbox{\rm and} \ \ \
\sum_{n \geq 1} \frac{N_{w,2}(n)}{2n(2n+1)(2n+2)},
$$
thus generalizing Corollary~1 in \cite{Sondow2} and a series for Euler's 
constant in \cite{Behrmann, Addison, Lint} (dated February 1967, 
August 1967, February 1968), respectively. Finally, 
we indicate some generalizations of our results, including an extension 
to base $B > 2$, and a method for giving alternate proofs without using 
the general lemma from \cite{AS}.

\section{A general lemma}

The first lemma in this section is taken from \cite{AS}; for completeness we 
recall the proof. We also give two classical results presented as lemmas,
together with a new result (Lemma~\ref{mainlemma}).

We start with some definitions.
Let $B \geq 2$ be an integer. Let $w$ be a word on the alphabet of digits
$\{0, 1, \cdots, B-1\}$ (that is, $w$ is a finite block of digits).
We denote by $N_{w,B}(n)$ the number of (possibly overlapping) occurrences
of $w$ in the $B$-ary expansion of an integer $n > 0$, and we set 
$N_{w,B}(0) = 0$.

Given $w$ as above, we denote by $|w|$ the length of the word $w$ (i.e., 
if $w = d_1 d_2 \cdots d_k$, then $|w| = k$). Denote by $w^j$ the 
concatenation of $j$ copies of the word $w$.

Given $w$ and $B$ as above, we denote by $v_B(w)$ the value of $w$ when $w$ 
is interpreted as the base $B$-expansion (possibly with leading $0$'s) of an
integer.

\begin{remark}
The occurrences of a given word in the $B$-ary expansion of the integer $n$ 
may overlap. For example, $N_{11,2}(7) = 2$.

If the word $w$ begins with $0$, but $v_B(w) \neq 0$, then in computing
$N_{w,B}(n)$ we assume that the $B$-ary expansion of $n$ starts with an 
arbitrarily long prefix of $0$'s. If $v_B(w) = 0$ we use the usual $B$-ary 
expansion of $n$ without leading zeros.
For example, $N_{011,2}(3) = 1$
(write $3$ in base $2$ as $0\cdots011$) and $N_{0,2}(2) = 1$.
\end{remark}

\begin{lemma}[\cite{AS}]\label{AS}
Fix an integer $B \geq 2$, and let $w$ be a non-empty word on the alphabet
$\{0, 1, \cdots, B-1\}$. If $f: {\mathbb N} \to {\mathbb C}$ is a function
with the property that $\sum_{n \geq 1} |f(n)| \log n < \infty$, then
$$
\sum_{n \geq 1} N_{w,B}(n) \left(f(n) - \sum_{0 \leq j < B} f(Bn+j)\right)
= \sum f(B^{|w|} n + v_B(w)),
$$
where the last summation is over $n \geq 1$ if $w=0^j$ for some $j \geq 1$,
and over $n \geq 0$ otherwise.
\end{lemma}

\proof (See \cite{AS}.) 
As $N_{w,B}(n) \leq \lfloor\frac{\log n}{\log B}\rfloor + 1$, all series
$\sum N_{w,B}(un+v) f(un+v)$, where $u$ and $v$ are nonnegative integers,
are absolutely convergent. Let $\ell$ be the last digit of $w$, and
let $g :=B^{|w|-1}$. Then
$$
\sum_{n \geq 0} N_{w,B}(n) f(Bn+\ell) =
\sum_{0 \leq k < g} \sum_{n \geq 0} N_{w,B}(gn+k) f(Bgn+Bk+\ell)
$$
and
$$
\sum_{n \geq 0} N_{w,B}(Bn+\ell) f(Bn+\ell) =
\sum_{0 \leq k < g} \sum_{n \geq 0} N_{w,B}(Bgn+Bk+\ell) f(Bgn+Bk+\ell).
$$
Now, if either $n \neq 0$ or $v_B(w) \neq 0$, then for $k = 0,1,\ldots, g-1$
we have
$$
N_{w,B}(Bgn+Bk+\ell) - N_{w,B}(gn+k) =
\left\{
\begin{array}{ll}
1, \ &\mbox{\rm if } k = \lfloor \frac{v_B(w)}{B} \rfloor; \\
& \\
0, \ &\mbox{\rm otherwise}.
\end{array}
\right.
$$
On the other hand, if $n=0$ and $v_B(w)=0$ (hence $\ell=0$), then the 
difference equals $0$ for every $k \in \{0, 1, \cdots, g-1\}$. Hence
$$
\begin{array}{lll}
\displaystyle\sum_{n \geq 0} N_{w,B}(Bn+\ell) f(Bn+\ell) -
\displaystyle\sum_{n \geq 0} N_{w,B}(n) f(Bn+\ell) &=& 
\displaystyle\sum f\left(Bgn + B \lfloor\frac{v_B(w)}{B}\rfloor + \ell \right)
\\ 
&=& \displaystyle\sum f(B^{|w|} n + v_B(w)) \ \ \ \ (*),
\end{array}
$$
the last two summations being over $n \geq 0$ if $w$ is not of the form $0^j$, 
and over $n \geq 1$ if $w = 0^j$ for some $j \geq 1$. We then write
$$
\begin{array}{lll}
\displaystyle\sum_{n \geq 0} N_{w,B}(n) f(n) &=&
\displaystyle\sum_{0 \leq j < B} \ \sum_{n \geq 0} N_{w,B}(Bn+j) f(Bn+j) \\
&=& \displaystyle\sum_{j \in [0,B)\setminus\{\ell\}} 
\ \displaystyle\sum_{n \geq 0} N_{w,B}(Bn+j) f(Bn+j)
+ \displaystyle\sum_{n \geq 0} N_{w,B}(Bn+\ell) f(Bn+\ell)
\end{array}
$$
which together with $(*)$ gives
$$
\sum_{n \geq 0} N_{w,B}(n) \left(
f(n) - \sum_{0 \leq j < B} f(Bn+j) \right)
= \sum f(B^{|w|}n+v_B(w)).
$$
Since $N_{w,B}(0) = 0$, the proof is complete. \endpf

     Now let $\Gamma$ be the usual gamma function, let
$\Psi := \Gamma'/\Gamma$ be the logarithmic derivative of the gamma
function, let $\zeta(s)$ be the Riemann zeta function, let
$\zeta(s,x) := \sum_{n \geq 0} (n+x)^{-s}$ be the Hurwitz zeta function, and
let $\gamma$ denote Euler's constant.

\begin{lemma}\label{Hurwitz}
If $a$ and $b$ are positive real numbers, then
$$
\sum_{n \geq 1} \left(\frac{1}{an} - \frac{1}{an+b}\right) =
\frac{1}{b} + \frac{\gamma+\Psi(b/a)}{a} .
$$
\end{lemma}

\proof
We write
$$
\begin{array}{lll}
\displaystyle\sum_{n \geq 1} \left(\frac{1}{an} - \frac{1}{an+b}\right) &=& 
\displaystyle\lim_{s \to 1_+} \sum_{n \geq 1}
\left(\frac{1}{(an)^s} - \frac{1}{(an+b)^s}\right)
= \frac{1}{a} \displaystyle\lim_{s \to 1_+} \sum_{n \geq 1}
\left(\frac{1}{n^s} - \frac{1}{(n+\frac{b}{a})^s}\right) \\
&=&
\displaystyle\frac{1}{a} \displaystyle\lim_{s \to 1_+} \left(\zeta(s) - 
\zeta\left(s,\frac{b}{a}\right) + \left(\frac{a}{b}\right)^s\right) \\
&=& \displaystyle\frac{1}{b} +
\frac{1}{a} \displaystyle\lim_{s \to 1_+} 
\left(\left(\zeta(s) - \frac{1}{s-1}\right)
- \left(\zeta\left(s,\frac{b}{a}\right) - \frac{1}{s-1}\right)\right) \\
&=&
\displaystyle\frac{1}{b} + 
\frac{1}{a}\left(\gamma + \frac{\Gamma'(b/a)}{\Gamma(b/a)}\right)
= \frac{1}{b} + \frac{\gamma + \Psi(b/a)}{a}
\end{array}
$$
(see for example \cite[p. 271]{WW}). \endpf

\begin{lemma}\label{Gamma}
For $x > 0$ we have
$$
\sum_{r \geq 1} \left(\frac{x}{r} - \log \left(1+\frac{x}{r}\right)\right)
= \log x + \gamma x + \log \Gamma(x).
$$
\end{lemma}

\proof Take the logarithm of the Weierstra{\ss} product for $1/\Gamma(x)$
(see, for example, \cite[Section~1.1]{ErMaObTr} or 
\cite[Section~12.1]{WW}). \endpf

The next lemma in this section is the last step before proving 
our theorems.

\begin{lemma}\label{mainlemma}
Let $a$ and $b$ be positive real numbers. Then
$$
\sum_{n \geq 1} \left(\frac{1}{an} - \log \frac{an+1}{an}\right)
= \log \Gamma\left(\frac{1}{a}\right) + \frac{\gamma}{a} - \log a
$$
and
$$
\sum_{n \geq 0} \left(\frac{1}{an+b} - \log \frac{an+b+1}{an+b}\right)
= \log \Gamma\left(\frac{b+1}{a}\right) - \log \Gamma\left(\frac{b}{a}\right) 
- \frac{\Psi(b/a)}{a}.
$$
\end{lemma}

\proof The proof is straightforward. The first formula follows directly from 
Lemma~\ref{Gamma}. To prove the second, write the $n$th term of the series for
$n \geq 1$ as the following sum of $n$th terms of three absolutely 
convergent series:
$$
\frac{1}{an+b} - \frac{1}{an} - \frac{b}{an} 
+ \log\left(1+\frac{b}{an}\right) + \frac{b+1}{an}
- \log\left(1+\frac{b+1}{an}\right);
$$
then use Lemmas~\ref{Hurwitz} and \ref{Gamma}.
\endpf

\section{Two theorems}

In this section we give two theorems that are consequences of 
Lemma~\ref{AS}, and that generalize results in \cite{Sondow2} 
and \cite{Behrmann, Addison, Lint}. 

\begin{theorem}\label{degree2}
Let $w$ be a non-empty word on the alphabet $\{0, 1\}$, and let 
$\Psi$ denote the logarithmic derivative of the Gamma function.

\begin{itemize}
\item[(a)] If $v_2(w) = 0$, then
$$
\sum_{n \geq 1} \frac{N_{w,2}(n)}{2n(2n+1)} =
\log \Gamma\left(\frac{1}{2^{|w|}}\right) + \frac{\gamma}{2^{|w|}}
- |w| \log 2.
$$

\item[(b)] If $v_2(w) \neq 0$, then
$$
\sum_{n \geq 1} \frac{N_{w,2}(n)}{2n(2n+1)} =
\log \Gamma\left(\frac{v_2(w)+1}{2^{|w|}}\right) - 
\log \Gamma\left(\frac{v_2(w)}{2^{|w|}}\right) -
\frac{1}{2^{|w|}}\Psi\left(\frac{v_2(w)}{2^{|w|}}\right).
$$
\end{itemize}
\end{theorem}

\proof
Let 
$$
A_n := \frac{1}{n} - \log \frac{n+1}{n}
$$
for $n \geq 1$. 
Noting that $A_n - A_{2n} - A_{2n+1} = \frac{1}{2n(2n+1)}$, the theorem
follows from Lemma~\ref{AS} with $B=2$, and $f(n) := A_n$ for $n \geq 1$, 
together with Lemma~\ref{mainlemma}. \endpf

\begin{example}\label{gammaplusminus}

Taking $w = 0$ and $w = 1$, and recalling that 
$\Gamma(1/2) = \sqrt \pi$ and $\Psi(1/2) = - \gamma -2 \log 2$ 
by Gau{\ss}'s theorem (see for example \cite[p. 19]{ErMaObTr} or 
\cite[p. 94]{Knuth}), we get
$$
\sum_{n \geq 1} \frac{N_{0,2}(n)}{2n(2n+1)} =
\frac{1}{2} \log \pi + \frac{\gamma}{2} - \log 2
$$
and
$$
\sum_{n \geq 1} \frac{s_2(n)}{2n(2n+1)} =
\sum_{n \geq 1} \frac{N_{1,2}(n)}{2n(2n+1)} = 
- \frac{1}{2} \log \pi + \frac{\gamma}{2} + \log 2.
$$
These equalities imply the formulas in the Introduction:
$$
\sum_{n \geq 1} \frac{N_{1,2}(n) \pm N_{0,2}(n)}{2n(2n+1)} = \gamma^{\pm}
$$
where (following the notations of \cite{Sondow2}) $\gamma^+ := \gamma$
and $\gamma^- := \log \frac{4}{\pi}$, which is Corollary~1
of \cite{Sondow2}.

\end{example}

\begin{remark}
The formulas in Theorem~\ref{degree2} are analogous to those in 
\cite[p.~25]{AS}.  The analogy becomes more striking if one uses 
Gau{\ss}'s theorem to write all expressions of the form $\Psi(x)$, 
with $x$ a rational number in $(0,1]$, using only trigonometric 
functions, logarithms, and Euler's constant. 
\end{remark}

\begin{theorem}\label{degree3}
Let $w$ be a non-empty word on the alphabet $\{0, 1\}$.

\begin{itemize}
\item[(a)] If $v_2(w) = 0$, then
$$
\sum_{n \geq 1} \frac{N_{w,2}(n)}{2n(2n+1)(2n+2)} =
\log \Gamma\left(\frac{1}{2^{|w|}}\right) + \frac{\gamma}{2^{|w|+1}}
- |w| \log 2 - \frac{1}{2^{|w|+1}}\Psi\left(\frac{1}{2^{|w|}}\right)
- \frac{1}{2}.
$$

\item[(b)] If $v_2(w) \neq 0$, then
$$
\begin{array}{ll}
\displaystyle\sum_{n \geq 1} \frac{N_{w,2}(n)}{2n(2n+1)(2n+2)} &=
\log \Gamma\displaystyle\left(\frac{v_2(w)+1}{2^{|w|}}\right) -
\log \Gamma\displaystyle\left(\frac{v_2(w)}{2^{|w|}}\right) \\
& \ \ \ \ \ - \displaystyle\frac{1}{2^{|w|+1}} 
\left(\Psi\left(\frac{v_2(w)}{2^{|w|}}\right)
+ \Psi\left(\frac{v_2(w)+1}{2^{|w|}}\right)\right).
\end{array}
$$
\end{itemize}
\end{theorem}

\proof 
Noting that $\frac{1}{2n(2n+1)} - \frac{1}{4} \cdot \frac{1}{n(n+1)} =
\frac{1}{2n(2n+1)(2n+2)}$, it suffices to use Theorem~\ref{degree2} 
and the following result, deduced from \cite[top of p. 26]{AS} in the 
case $B=2$.

{\it
\begin{itemize}
\item[(a)]
If $v_2(w) = 0$, then
$$
\sum_{n \geq 1} \frac{N_{w,2}(n)}{n(n+1)} =
\frac{1}{2^{|w|-1}}
\left(
\Psi\left(\frac{1}{2^{|w|}}\right)
+ \gamma + 2^{|w|}
\right).
$$

\item[(b)] If $v_2(w) \neq 0$, then
$$
\sum_{n \geq 1} \frac{N_{w,2}(n)}{n(n+1)} =
\frac{1}{2^{|w|-1}}
\left(
\Psi\left(\frac{v_2(w)+1}{2^{|w|}}\right)
- \Psi\left(\frac{v_2(w)}{2^{|w|}}\right)
\right).
$$
\end{itemize}
}
\endpf

\begin{example}

Taking $w=0$ and $w=1$, we get
$$
\sum_{n \geq 1} \frac{N_{0,2}(n)}{2n(2n +1)(2n+2)} =
\frac{1}{2} \log \pi + \frac{\gamma}{2} - \frac{1}{2} \log 2 - \frac{1}{2}
$$
and
$$
\sum_{n \geq 1} \frac{s_2(n)}{2n(2n +1)(2n+2)} =
\sum_{n \geq 1} \frac{N_{1,2}(n)}{2n(2n +1)(2n+2)} = 
- \frac{1}{2} \log \pi + \frac{\gamma}{2} + \frac{1}{2} \log 2.
$$
Hence
$$
\sum_{n \geq 1} \frac{N_{1,2}(n) \pm N_{0,2}(n)}{2n(2n +1)(2n+2)} = \delta^{\pm}
$$
where $\delta^+ := \gamma -\frac{1}{2}$ and 
$\delta^- := \frac{1}{2} - \log \frac{\pi}{2}$, which are respectively 
a formula given in \cite{Behrmann, Addison, Lint} and a seemingly 
new companion formula.

\end{example}

\begin{remark}
As mentioned, all expressions of the form $\Psi(x)$, with
$x$ a rational number in $(0,1]$, can be written using only trigonometric
functions, logarithms, and Euler's constant.
\end{remark}

\section{Generalizations}\label{generalizations}

Several extensions or generalizations of our results are possible.
We give some of them in this section.

\subsection{Variation on $A_n$}

Instead of applying Lemma~\ref{AS} with $f(n) = A_n = \frac{1}{n} - \log \frac{n+1}{n}$ for $n \geq 1$,
we could replace $A_n$ with 
$$
A^{(k)}_n := \frac{1}{n+k} - \log \frac{n+1}{n}
$$
for $n \geq 1$, where $k$ is a nonnegative integer. Defining the (rational) 
function $Q^{(k)}$ by
$$
Q^{(k)}(n) := A^{(k)}_n - A^{(k)}_{2n} - A^{(k)}_{2n+1}
$$
and noting that summing $\sum_{n \geq 1} A^{(k)}_{an+b}$ boils down to 
summing $\sum_{n \geq 1} \left(\frac{1}{an+b+k} - \frac{1}{an+b}\right)$, 
which as in the proof of Lemma~\ref{Hurwitz} involves the Hurwitz zeta 
function, we obtain explicit formulas for the sum of the series 
$\sum_{n \geq 1} N_{w,2}(n) Q^{(k)}(n)$.

\subsection{Extension to base $B > 2$}\label{gal}

Lemma~\ref{AS} has been used above only for base $B=2$. There are
applications to other bases in \cite{AS}. We also note that the relation
among the $A_n$'s,
$$
A_n = \frac{1}{2n(2n+1)} + A_{2n} + A_{2n+1}
$$
for $n \geq 1$, can be generalized to base $B$. Namely,
$$
A_n = Q(n,B) + R(n,B)
$$
where
$$
Q(n,B) := \frac{1}{Bn(Bn+1)} + \frac{2}{Bn(Bn+2)} + 
\cdots + \frac{B-1}{Bn(Bn+B-1)}
$$
and
$$
R(n,B) := A_{Bn} + A_{Bn+1} + \cdots + A_{Bn+B-1}.
$$
This allows us to use Lemmas~\ref{AS} and \ref{mainlemma} to sum, for 
example, the series
$$
\sum_{n \geq 1} N_{w,3}(n) \frac{9n+4}{3n(3n+1)(3n+2)},
$$
since
$$
Q(n,3) = A_n - A_{3n} - A_{3n+1} - A_{3n+2} = \frac{9n+4}{3n(3n+1)(3n+2)}.
$$

\subsection{Weighted $A_n$'s}

In this section we consider a weighted form of the $A_n$'s.  
First we need to study a relation between sequences of real numbers.

\begin{lemma}\label{transform}
Let $(r_n)_{n \geq 1}$ and $(R_i)_{i \geq 1}$ be sequences of real numbers.
Set $r_0 := 0$ and $R_0 := 0$. Then the following two properties are equivalent:
\begin{itemize}
\item[(1)] for $i \geq 1$
$$
R_i = \sum_{k \geq 0} r_{\lfloor \frac{i}{2^k} \rfloor} =
r_i + r_{\lfloor\frac{i}{2}\rfloor} + r_{\lfloor\frac{i}{4}\rfloor} 
+ \cdots
$$
(note that this is actually a finite sum);

\item[(2)] for $n \geq 1$
$$
r_n = R_n - R_{\lfloor\frac{n}{2}\rfloor}.
$$
\end{itemize}
\end{lemma}

\proof 
The implication (1) $\Rightarrow$ (2) is easily seen by considering the 
cases $n$ even and $n$ odd. Likewise, for (2) $\Rightarrow$ (1) take 
$i$ even and $i$ odd.  \endpf

\begin{remark}
See \cite[Theorem~9]{ASreg2} for more about this relation.
\end{remark}

\begin{theorem}\label{weight}
Assume that $r_1, r_2, \ldots$ and $R_1, R_2, \ldots$ are real 
numbers related as in Lemma~\ref{transform}. Then the series 
$\sum |r_n| n^{-2}$ converges if and only if the series
$\sum |R_i| i^{-2}$ converges, and in this case we have
$$
S := \sum_{n \geq 1} r_n \left(\frac{1}{n} - \log \frac{n+1}{n}\right)
= \sum_{i \geq 1} \frac{R_i}{2i(2i+1)}
$$
\end{theorem}

\proof  First note that if the series $\sum |R_i|i^{-2}$ converges, 
then so does the series $\sum |r_n|n^{-2}$: use the expression for 
$r_n$ in terms of the $R_i$'s in Lemma~\ref{transform}. 
Now suppose that the series 
$\sum |r_n| n^{-2}$ converges. As before, let 
$A_n := \frac{1}{n} - \log \frac{n+1}{n}$. 
Then $0 < A_n < \frac{1}{n} - \frac{1}{n+1}$.
This implies that the series 
$S :=\sum r_n \left(\frac{1}{n} - \log \frac{n+1}{n}\right)$
is absolutely convergent. Now
$$
A_n = \frac{1}{2n(2n+1)} + A_{2n} + A_{2n+1}
$$ 
implies 
$$
A_n = \frac{1}{2n(2n+1)} + \frac{1}{4n(4n+1)} + \frac{1}{(4n+2)(4n+3)} 
+ A_{4n} + A_{4n+1} + A_{4n+2} + A_{4n+3}.
$$
Hence, repeating $K$ times,
$$
A_n = \sum_{1 \leq k \leq K} \sum_{0 \leq m < 2^{k-1}}
\frac{1}{(2^k n + 2m)(2^k n + 2m + 1)} 
+ \sum_{0 \leq q < 2^K} A_{2^K n + q}.
$$
Using the bounds $0 < A_n < \frac{1}{n} - \frac{1}{n+1}$ and
telescoping, the last sum is less than $2^{-K}$. Letting $K$
tend to infinity, we obtain the (rapidly convergent) series
$$
A_n = \sum_{k \geq 1} \ \sum_{0 \leq m < 2^{k-1}}
\frac{1}{(2^k n + 2m)(2^k n + 2m + 1)}.
$$
Substituting into the sum defining $S$ yields the double series
$$
S = \sum_{n \geq 1} \ \sum_{k \geq 1} \ \sum_{0 \leq m < 2^{k-1}}
\frac{r_n}{(2^k n + 2m)(2^k n + 2m + 1)},
$$
which converges absolutely. Thus we may collect terms with
the same denominator, and we arrive at the series
$$
S = \sum_{i \geq 1} \frac{R'_i}{2i(2i+1)},
$$
where
$$
R'_{i} := \sum_{n \in {\cal E}_{i}} r_n
$$
with ${\cal E}_{i} := \{n \geq 1, \ \exists k \geq 1, \
\exists m \in [0, 2^{k-1}), \ 2^{k-1}n + m = i \}$.
On the one hand, this proves that the series 
$\sum \frac{R'_i}{2i(2i+1)}$ is absolutely convergent 
(hence the series $\sum |R'_i| i^{-2}$ is convergent).
On the other hand, $R'_i$ can also be written as
$$
R'_{i} := \sum_{1 \leq k \leq \frac{\log i}{\log 2} + 1} 
r_{\lfloor\frac{i}{2^{k-1}}\rfloor} = 
\sum_{k \geq 0} r_{\lfloor\frac{i}{2^k}\rfloor}
$$
(recall that we have set $r_0 := 0$). Thus $R'_i = R_i$ 
by the hypothesis, and the proof is complete. \endpf

\begin{example}\label{fouroverpi}
Theorem~\ref{weight} yields in particular the series for $\gamma$ and 
$\log \frac{4}{\pi}$ in the Introduction, in Example~\ref{gammaplusminus},
and in \cite[Corollary 1]{Sondow2}. Namely,

If $r_1 = r_2 = \cdots = 1$, then the series defining $S$ sums to 
$\gamma$ from Lemma~\ref{mainlemma}, and the formula defining $R'_i=R_{i}$ 
reduces to $R_{i} = \lfloor \frac{\log 2i}{\log 2} \rfloor = 
N_{1,2}(i) + N_{0,2}(i)$.

If $r_n = (-1)^{n-1}$, then $S = \log \frac{4}{\pi}$ 
(see \cite{Sondow1} or decompose $S$ into 
$\sum$ (odd terms) 
$- \sum$ (even terms) and apply Lemma~\ref{mainlemma}), and the 
formula defining $R'_i=R_{i}$ implies $R_{i} = N_{1,2}(i) - N_{0,2}(i)$.
To see this equality, first note that if it holds for $i \geq 1$, then 
using Lemma 5 and looking at the cases $n$ odd and $n$ even,
$$
r_n = R_n - R_{\left\lfloor \frac{n}{2} \right\rfloor} = (-1)^{n-1}
$$
for $n \geq 1$ (compare \cite[Lemma~2]{Sondow2}). Now recall that 
properties (1) and (2) in Lemma~\ref{transform} are equivalent.

\end{example}

\begin{remark}\label{newproof}

Example~\ref{fouroverpi} shows that it is possible to deduce the formula
$$
\sum_{n \geq 1} \frac{N_{1,2}(n)-N_{0,2}(n)}{2n(2n+1)} = \log \frac{4}{\pi}
$$
 from Theorem~\ref{weight} and Lemma~\ref{mainlemma} without using
Lemma~\ref{AS}: this yields a proof of the formula that is different 
from those in \cite{Sondow2} and Example~\ref{gammaplusminus}. Similar 
reasoning applies for any ultimately periodic sequence $(r_n)_{n \geq 1}$. 
In particular, it is not hard to see that the relations giving $r_{2n}$ and 
$r_{2n+1}$ in terms of the $R_i$'s imply that the sequence $(r_n)_{n \geq 1}$ 
is periodic whenever $R_i = N_{w,2}(i)$ for some fixed $w$ and for every 
$i \geq 1$. Hence Theorem~\ref{degree2} can be deduced from Theorem~\ref{weight} 
and Lemma~\ref{mainlemma} (along with the method for decomposing series employed
in Example~\ref{fouroverpi}), without using Lemma~\ref{AS}. In the same vein, 
the generalization in Section~\ref{gal} can be proved using a
generalization of Theorem~\ref{weight} to base $B$ together with
Lemma~\ref{mainlemma}.

\end{remark}

\section{Future directions}

Lemma~\ref{AS} is the main tool for summing series in \cite{AS} and in the 
present paper. It might be possible to use the lemma to obtain the base $B$ 
accelerated series for Euler's constant in \cite[Theorem 2]{Sondow2}, and 
to sum more general series with $N_{w,B}(n)$. On the other hand, it might 
also be possible to extend the results of \cite{AS} and the present paper, 
and sum series where $(N_{w,B}(n))_{n \geq 1}$ is replaced by a more
general integer sequence $(a_n)_{n \geq 1}$, using the decomposition in 
\cite{MM} of a sequence $(a_n)_{n \geq 1}$ into a (possibly infinite) 
linear combination of block-counting sequences $(N_{w,B}(n))_{n \geq 1}$ 
(see also \cite{ASreg}). Of course, since this may replace a series with 
an infinite sum, for the method to work the new series must be summable 
in closed form.

\end{document}